\documentclass{article}
\pdfoutput=1
\usepackage{amsmath}
\usepackage{amssymb}
\usepackage{hyperref}
\usepackage{tikzexternal}
\tikzsetexternalprefix{figures/}
\tikzsetfigurename{main-}
\providecommand{\tikzifexternalizing}[2]{#2}
\newtheorem{theorem}{Theorem}[section]
\newtheorem{lemma}[theorem]{Lemma}
\newtheorem{problem}[theorem]{Problem}
\newtheorem{proposition}[theorem]{Proposition}
\begin{document}

\title{Benjamini-Schramm convergence and the distribution of chromatic roots for
sparse graphs}
\author{Mikl\'{o}s Ab\'{e}rt and Tam\'{a}s Hubai}
\maketitle

\begin{abstract}
We define the chromatic measure of a finite simple graph as the uniform
distribution on its chromatic roots. We show that for a Benjamini-Schramm
convergent sequence of finite graphs, the chromatic measures converge in
holomorphic moments.

As a corollary, for a convergent sequence of finite graphs, we prove that the
normalized log of the chromatic polynomial converges to an analytic function
outside a bounded disc. This generalizes a recent result of Borgs, Chayes,
Kahn and Lov\'{a}sz, who proved convergence at large enough positive integers
and answers a question of Borgs.

Our methods also lead to explicit estimates on the number of proper colorings
of graphs with large girth.

\end{abstract}

\section{Introduction}

Let $G$ be a finite undirected graph without multiple edges and loops. A map
$f:V(G)\rightarrow\left\{1,\ldots,q\right\}$ is a \emph{proper coloring}
if for all edges $(x,y)\in E(G)$ we have $f(x)\neq f(y)$. For a positive
integer $q$ let $\mathrm{ch}_{G}(q)$ denote the number of proper colorings of
$G$ with $q$ colors. Then $\mathrm{ch}_{G}$ is a polynomial in $q$, called the
\emph{chromatic polynomial} of $G$. The complex roots of $\mathrm{ch}_{G}$ are
called the \emph{chromatic roots} of $G$.

The study of chromatic polynomials and its roots has been initiated by
Birkhoff. Since then, the subject has gotten considerable interest, partially because
of its connection to statistical mechanics. In particular, the chromatic roots
control the behaviour of the antiferromagnetic Potts model at zero
temperature. For a survey on the subject see \cite{sokalsurvey}.

For a finite graph $G$, a finite rooted graph $\alpha$ and a positive integer
$R$, let $\mathbf{P}(G,\alpha,R)$ denote the probability that the $R$-ball
centered at a uniform random vertex of $G$ is isomorphic to $\alpha$. We say
that a graph sequence $(G_{n})$ of bounded degree is \emph{Benjamini-Schramm
convergent} if for all finite rooted graphs $\alpha$ and $R>0$, the
probabilities $\mathbf{P}(G_{n},\alpha,R)$ converge (see \cite{bensch}). This means
that one can not statistically distinguish $G_{n}$ and $G_{n^{\prime}}$ for
large $n$ and $n^{\prime}$ by sampling them from a random vertex with a fixed
radius of sight. An example (that is regularly used in statistical physics) is
to approximate the infinite lattice $\mathbb{Z}^{d}$ by bricks with all the
side lengths tending to infinity. More generally, amenable vertex transitive
graphs can be obtained as the Benjamini-Schramm limits of their F\o lner sequences.

For a simple graph $G$ let $\mu_{G}$, the \emph{chromatic measure of} $G$
denote the uniform distribution on its chromatic roots. By a theorem of Sokal
\cite{sokal}, $\mu_{G}$ is supported on the open disc of radius $Cd$ around $0$,
denoted by
\[
D=B(0,Cd)
\]
where $d$ is the maximal degree of $G$ and $C<8$ is an absolute constant.

\begin{theorem}
\label{moments}Let $(G_{n})$ be a Benjamini-Schramm convergent graph sequence
of absolute degree bound $d$, and $\widetilde{D}$ an open neighborhood of
the closed disc $\overline{D}$. Then for every holomorphic function
$f:\widetilde{D}\rightarrow\mathbb{C}$, the sequence
\[
{\displaystyle\int\limits_{D}}
f(z)d\mu_{G_{n}}(z)
\]
converges.
\end{theorem}

Let $\ln$ denote the principal branch of the complex logarithm function. For a
simple graph $G$ and $z\in\mathbb{C}$ let
\[
\mathrm{t}_{G}(z)=\frac{\ln \mathrm{ch}_{G}(z)}{\left\vert V(G)\right\vert}
\]
where this is well-defined. In statistical mechanics, $\mathrm{t}_{G}(z)$ is
called the \emph{entropy per vertex} or the \emph{free energy} at $z$. In
their recent paper \cite{lovasz}, Borgs, Chayes, Kahn and Lov\'{a}sz proved
that if $(G_{n})$ is a Benjamini-Schramm convergent graph sequence of absolute
degree bound $d$, then $\mathrm{t}_{G_{n}}(q)$ converges for every positive
integer $q>2d$. Theorem \ref{moments} yields the following.

\begin{theorem}
\label{holom}Let $(G_{n})$ be a Benjamini-Schramm convergent graph sequence of
absolute degree bound $d$ with $\left\vert V(G_{n})\right\vert \rightarrow\infty
$. Then $\mathrm{t}_{G_{n}}(z)$ converges to a real analytic function on
$\mathbb{C}\setminus\overline{D}$.
\end{theorem}

In particular, $\mathrm{t}_{G_{n}}(z)$ converges for all $z\in\mathbb{C}
\setminus\overline{D}$. Theorem \ref{holom} answers a question of Borgs
\cite[Problem 2]{borgs} who asked under which circumstances the entropy per
vertex has a limit and whether this limit is analytic in $1/z$. Note that for
an amenable quasi-transitive graph and its F\o lner sequences, this was shown
to hold in \cite{ize}.

To prove Theorem \ref{moments} we show that for a finite graph $G$ and for
every $k$, the number
\[
p_{k}(G)=\left\vert V(G)\right\vert
{\displaystyle\int\limits_{D}}
z^{k}d\mu_{G}(z)
\]
can be expressed as a fixed linear combination of $\hom(H,G)$ where the $H$
are connected finite graphs and $\hom(H,G)$ denotes the number of graph
homomorphisms from $H$ to $G$. Since a sequence of graphs $G_{n}$ of bounded
degree is Benjamini-Schramm convergent if and only if
\[
\frac{\hom(H,G_{n})}{\left\vert V(G_{n})\right\vert}
\]
converges for all connected graphs $H$. This gives convergence of all the
holomorphic moments of $\mu_{G_{n}}$, and this is equivalent to Theorem
\ref{moments}. For instance, for the fourth moment we get

\noindent$p_{4}(G)=\input{includes/p4e.inc}$.\bigskip

One could speculate that assuming Benjamini-Schramm convergence of $G_{n}$,
maybe the complex measures $\mu_{G_{n}}$ themselves will weakly converge. That
is, Theorem \ref{moments} would hold for any continuous real function on $D$
or, equivalently, convergence would hold in all the moments
\[
{\displaystyle\int\limits_{D}}
z^{k}\overline{z}^{j}d\mu_{G_{n}}(z)\text{.}
\]

However, this is not true in general, as the following easy counterexample
shows. Let $P_{n}$ denote the path of length $n$ and let $C_{n}$ denote the
cycle of length $n$. Then $P_{n}$ and $C_{n}$ converge to the same object, the
infinite rooted path, while we have
\[
\mathrm{ch}_{P_{n}}(z)=z(z-1)^{n-1}\text{ and }\mathrm{ch}_{C_{n}}(z)=(z-1)^{n}+(-1)^{n}
(z-1)\text{.}
\]
Thus, the weak limit of $\mu_{P_{n}}$ is the Dirac measure on $1$ and the weak
limit of $\mu_{C_{n}}$ is the normalized Lebesgue measure on the unit circle
centered at $1$.

Still, using Theorem \ref{moments}, we are able to prove the weak convergence
of $\mu_{G_{n}}$ for some natural sequences of graphs. For example, let
$T_{n}=C_{4}\times P_{n}$ denote the $4\times n$ tube, i.e.\ the cartesian
product of the 4-cycle with a path on $n$ vertices. $T_{n}$ is a $4$-regular
graph except at the ends of the tube.

\begin{proposition}
\label{tube}The chromatic measures $\mu_{T_{n}}$ weakly converge.
\end{proposition}

The proof is as follows: as Salas and Sokal \cite{salas} showed, the pointwise
limit $X$ of supports of $\mu_{T_{n}}$ is part of a particular algebraic
curve, so any subsequential weak limit is supported on $X$. The complement of
$X$ is connected, so by Mergelyan's theorem \cite{merg}, every continuous
function on $X$ can be uniformly approximated by polynomials. Using Theorem
\ref{moments} this yields weak convergence of $\mu_{T_{n}}$. See Section
\ref{convsection} for details.

In this case one can use the so-called transfer matrix method to control the
support of the chromatic measures (see \cite{salas} for various models related
to the square lattice). In general, even for models of the square lattice, the
complement of the limiting set may not be connected, and hence one can not
invoke Mergelyan's theorem. It is expected, however, that for any model where
the transfer matrix method can be used, the chromatic measures do converge weakly.

Another naturally interesting case is when the girth of $G$ (the minimal size
of a cycle) is large. One can show that
\[
{\displaystyle\int\limits_{D}}
z^{k}d\mu(z)=\frac{\left\vert E(G)\right\vert}{\left\vert V(G)\right\vert
}\text{ \ (}1\leq k\leq\mathrm{girth}(G)-2\text{)}
\]
that is, the moments are all the same until the girth is reached (see Lemma
\ref{girthmoments}). This implies that for a sequence of $d$-regular graphs
$G_{n}$ with girth tending to infinity, the limit of the free entropy
\[
\lim_{n\rightarrow\infty}\mathrm{t}_{G_{n}}(z)=\ln q+\frac{d}{2}\ln
(1-\frac{1}{q})
\]
for $q>Cd$. This is one of the main results in \cite{gamarnik}. Note that
their proof works for $q>d+1$, while our approach only works for $q>Cd$. The
advantage of our approach is that it gives an explicit estimate on the number
of proper colorings of large girth graphs.

\begin{theorem}
\label{girth}Let $G$ be a finite graph of girth $g$ and maximal degree $d$.
Then for all $q>Cd$ we have
\[
\left\vert \frac{\ln \mathrm{ch}_{G}(q)}{\left\vert V(G)\right\vert}-\left(\ln
q+\frac{\left\vert E(G)\right\vert}{\left\vert V(G)\right\vert}\ln
(1-\frac{1}{q})\right)\right\vert \leq2\frac{(Cd/q)^{g-1}}{1-Cd/q}\text{.}
\]
\end{theorem}

When $G$ is $d$-regular with asymptotically maximal girth, i.e.\ %
$g=c\ln\left\vert V(G)\right\vert$, this yields
\[
\left\vert \frac{\ln \mathrm{ch}_{G}(q)}{\left\vert V(G)\right\vert}-\left(\ln
q+\frac{d}{2}\ln(1-\frac{1}{q})\right)\right\vert
\leq O(\left\vert V(G)\right\vert ^{-c^{\prime}})
\]
for some explicit constant $c^{\prime}>0$. Counting the number of proper
colorings of random $d$-regular graphs have been considered in \cite{gamarnik}.
These graphs do not have logarithmic girth, but they have few shorter cycles,
so one can obtain a similar result for them.

Here we wish to raise attention to an interesting phenomenon, of which we only
have some computational evidence. We have computed the chromatic measures of
several $3$-regular large girth graphs and surprisingly, it looks like one
may also get weak convergence of chromatic measures.

\begin{figure}[h]
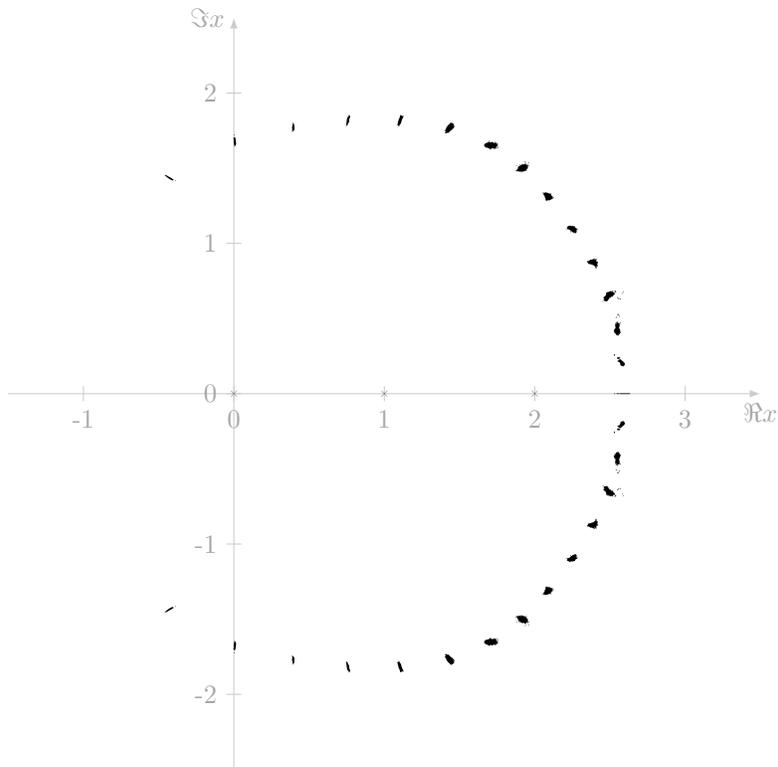

\begin{center}
\tikzsetnextfilename{lgirth}
\begin{tikzpicture}[scale=2]
\draw[thin,color=black!20,-latex]
    (-1.5,0) -- (3.5,0)
    node[below,color=black!35] {$\Re x$};
\draw[thin,color=black!20,-latex]
    (0,-2.5) -- (0,2.5)
    node[left,color=black!35] {$\Im x$};
\foreach \i in {-1,...,3}
    \draw[thin,color=black!20]
    (\i,0.05) -- (\i,-0.05)
    node[below,color=black!35] {\i};
\foreach \i in {-2,...,2}
    \draw[thin,color=black!20]
    (0.05,\i) -- (-0.05,\i)
    node[left,color=black!35] {\i};
\input{includes/crosses.inc}
\input{includes/circaprx.inc}
\end{tikzpicture}
\end{center}
\caption{Chromatic roots of the 30368 cubic graphs of size 32 and girth 7}
\label{largegirthpicture}
\end{figure}

\begin{problem}
Let $G_{n}$ be a sequence of $d$-regular graphs with girth tending to
infinity. Does $\mu_{G_{n}}$ weakly converge?
\end{problem}

This would be interesting because one could consider the limit as the
`chromatic measure of the $d$-regular infinite tree'. Observe that any
subsequential weak limit $\mu$ of $\mu_{G_{n}}$ satisfies
\[
{\displaystyle\int\limits_{D}}
z^{k}d\mu(z)=\frac{d}{2}\text{ \ (}k\geq1\text{)}
\]
that is, the holomorphic moments of $\mu$ are independent of $k$. While Figure
\ref{largegirthpicture} looks very promising, one misleading aspect of it may
be that $3$-regular graphs having $32$ vertices and (maximal possible) girth
$7$ may exhibit structural restrictions that are much stronger than just large
girth. \bigskip

It would be interesting to generalize the results of this paper to the Tutte
polynomial $T_{G}(x,y)$. This two-variable polynomial encodes a lot of
interesting invariants of $G$. For instance, $T_{G}(z,0)=\mathrm{ch}_{G}(z)$,
$T_{G}(2,1)$ counts the number of forests, $T_{G}(1,1)$ is the number of
spanning trees and $T_{G}(1,2)$ is the number of connected spanning subgraphs.
By a result of Lyons \cite{lyons}, we know that
\[
\frac{\log T_{G_{n}}(1,1)}{\left\vert V(G_{n})\right\vert}
\]
converges for a Benjamini-Schramm convergent sequence of graphs $G_{n}$ of
bounded degree. Also, in \cite{lovasz} it is shown that the same holds at the
places $(q,y)$ where $0\leq y<1$ and $q$ is large enough in terms of the
maximal degree. It would be interesting to see whether this also holds at
other places. The places $(2,1)$ and $(1,2)$ would be good test points as they
have a natural combinatorial meaning. Also, it is not clear whether Theorem
\ref{moments} holds for $p_{G}(z)=T_{G}(z,y_{0})$ for all fixed $y_{0}$. Note
that even for the chromatic polynomial, in general, the above log convergence
will not hold, for instance at $(2,0)$, because cycles of even and odd length
converge to the same limit, but even cycles have a proper $2$-coloring, while
odd cycles do not. This may not be so surprising, since $T_{G}(2,0)\leq
2^{c(G)}$ where $c(G)$ is the number of components of $G$. So for a nontrivial
graph sequence $G_{n}$, $T_{G_{n}}(2,0)$ is subexponential in $\left\vert
V(G_{n})\right\vert $, which points to the proximity of roots of $T_{G_{n}}$. To
apply Theorem \ref{moments} in its present form, one needs that some small
neighbourhood of the place is sparse in terms of roots. \bigskip

\noindent\textbf{Remark.} Note that recently Csikv\'{a}ri and Frenkel
\cite{csf} generalized Theorem \ref{holom} to a large class of graph
polynomials, including the Tutte polynomial.  In particular, they show that
convergence holds for the normalized log of $T(x,y)$ where $x,y$ have large
enough absolute value.\\

\noindent\textbf{Acknowledgements.} The authors thank L\'{a}szl\'{o}
Lov\'{a}sz who introduced them to \cite{lovasz}, G\'{a}bor Hal\'{a}sz who
raised their attention to the relevance of Mergelyan's theorem and
Lewis Bowen who pointed out some inconsistencies of notation in a previous
version of this paper. This work was partially supported by the MTA R\'{e}nyi
``Lend\"{u}let'' Groups and Graphs Research Group.

\section{Preliminaries}

For a simple graph $H$ on $n$ vertices let the number of legal colorings of
$H\ $with $q$ colors be denoted by $\mathrm{ch}_{H}(q)$. Then for any edge $e$
of $H$ the following recursion holds:
\[
\mathrm{ch}_{H}(q)=\mathrm{ch}_{H\setminus e}(q)-\mathrm{ch}_{H/e}(q)
\]
where $H\setminus e$ is obtained from $H$ by deleting $e$ and $H/e$ is
obtained by gluing the endpoints of $e$ and erasing multiple edges and loops.
This implies that $\mathrm{ch}_{H}$ is a polynomial of degree $n$ in $q$ with
integer coefficients, called the \emph{chromatic polynomial} of $G$ and that
the above recursion holds for the polynomials themselves. It also follows that
the constant coefficient of $\mathrm{ch}_{H}$ is zero and its main coefficient
is $1$. So, we can write
\begin{align*}
\mathrm{ch}_{H}(z)&=z^{n}-e_{1}(H)z^{n-1}+\ldots+(-1)^{k}e_{k}
(H)z^{n-k}+\ldots+(-1)^{n-1}e_{n-1}(H)z=\\
&={\displaystyle\prod\limits_{i=1}^{n}}
(z-\lambda_{i}(H))\text{.}
\end{align*}
The $e_{k}(H)$ are called the \emph{chromatic coefficients} of $H$ and
$\lambda_{i}(H)$ are its \emph{chromatic roots}. For $k\geq0$ let
\[
p_{k}(H)=
{\displaystyle\sum\limits_{i=1}^{n}}
\lambda_{i}^{k}(H)\text{.}
\]
The Newton identities establish connections between the roots and coefficients
of a polynomial. In this paper we will use the following version:
\[
p_{k}=(-1)^{k-1}ke_{k}+
{\displaystyle\sum\limits_{i=1}^{k-1}}
(-1)^{k-i-1}p_{i}e_{k-i}\text{.}
\]

Let $H,G$ be simple graphs. A map $f:V(H)\rightarrow V(G)$ is a
\emph{homomorphism} if for all edges $(x,y)\in E(H)$ we have $(f(x),f(y))\in
E(G)$. We denote the number of homomorphisms from $H$ to $G$ by $\hom(H,G)$.
The quantity $\hom(H,G)$ is nice to work with, mainly because of the following property.

\begin{lemma}
\label{szorzat}Let $H$ be the disjoint union of $H_{1}$ and $H_{2}$. Then
\[
\hom(H,G)=\hom(H_{1},G)\hom(H_{2},G)
\]
for all simple graphs $G$.
\end{lemma}

We leave the proof to the reader.\bigskip

For a random rooted graph $G$, a finite rooted graph $\alpha$ and a positive integer
$R$, let $\mathbf{P}(G,\alpha,R)$ denote the probability that the $R$-ball
centered at the root of $G$ is isomorphic to $\alpha$. Analogously, for an unrooted finite
graph $G$, let $\mathbf{P}(G,\alpha,R$) denote the probability that the $R$-ball
centered at a uniform random vertex of $G$ is isomorphic to $\alpha$. A graph
sequence $G_{n}$ has bounded degree if there is an absolute upper bound on
the degrees of vertices of $G_{n}$.

A graph sequence $(G_{n})$ of bounded degree is \emph{Benjamini-Schramm
convergent} if for all finite rooted graphs $\alpha$ and $R>0$ the
probabilities $\mathbf{P}(G_{n},\alpha,R)$ converge.

The limit of a Benjamini-Schramm convergent sequence of graphs is the random
rooted graph $G$ satisfying
\[
\mathbf{P}(G,\alpha,R)=\lim_{n\rightarrow\infty}\mathbf{P}(G_{n},\alpha,R)
\]
for all $R>0$ and $\alpha$. It is easy to see that $G$ is well defined. In the
most transparent case, $G$ is just one graph, which then has to be vertex
transitive. For instance, the $d$ dimensional lattice
\[
\mathbb{Z}^{d}=\lim_{n\rightarrow\infty}(\mathbb{Z}/n\mathbb{Z)}^{d}
=\lim_{n\rightarrow\infty}B_{n}^{d}
\]
where $(\mathbb{Z}/n\mathbb{Z)}^{d}$ is the $d$ dimensional torus and
$B_{n}^{d}$ is the box of side length $n$ in $\mathbb{Z}^{d}$. The same way,
one can obtain any connected vertex transitive amenable graph as a limit. Let
$G$ be a connected vertex transitive graph of bounded degree. A sequence of
connected subgraphs $F_{n}$ of $G$ is a \emph{F\o lner sequence}, if
\[
\lim_{n\rightarrow\infty}\frac{\left\vert \partial F_{n}\right\vert
}{\left\vert V(F_{n})\right\vert}=0
\]
where $\partial F_{n}$ denotes the external vertex boundary of $F_{n}$. Note
that $G$ is amenable if and only if it has a F\o lner sequence. It is easy to
show that any connected vertex transitive amenable graph is the
Benjamini-Schramm limit of its F\o lner sequences.

Let us consider now the $d$-regular tree $T_{d}$, which is in many senses the
farthest possible from being amenable. One can obtain $T_{d}$ as the limit of
finite graphs, but it is worth to point out that $T_{d}$ can \emph{not} be
obtained as a limit of finite trees. Indeed, the expected degrees of finite
trees are approximately $2$ and this passes on to their limits. It is a good
exercise to understand what the limit of the balls in $T_{d}$ is (it is a
fixed tree where the root is random). The right way to approximate $T_{d}$ in
Benjamini-Schramm convergence is to take finite $d$-regular graphs $G_{n}$
with girth tending to infinity.

Benjamini-Schramm convergence can also be expressed in terms of graph
homomorphisms using the following lemma (see \cite{lovaszbook}, Proposition 5.6).

\begin{lemma}
\label{convhom}Let $G_{n}$ be a graph sequence of bounded degree. Then $G_{n}$
is Benjamini-Schramm convergent if and only if for every finite connected
graph $H$, the limit
\[
\lim_{n\rightarrow\infty}\frac{\hom(H,G_{n})}{\left\vert V(G_{n})\right\vert}
\]
exists.
\end{lemma}

Note that one needs connectedness in Lemma \ref{convhom}, as $\hom(H,G)$ may
be the order of $\left\vert V(G)\right\vert ^{c}$ where $c$ is the number of
components in $H$.

\section{Expressing moments from homomorphisms}

In this section we give an explicit formula for the holomorpic moments of the
chromatic measure in terms of graphs homomorphisms. \bigskip

For a finite, simple graph $G$ let $\mathcal{P(}G\mathcal{)}$ be defined as
the set of partitions of $V(G)$ where no edge of $G$ connects two vertices in
the same class. A partition $P\in\mathcal{P(}G\mathcal{)}$ can be considered
as a surjective homomorphism from $G$ to the simple graph $G/P$ obtained by
contracting each class of $P$ and erasing multiple edges. For simple graphs
$G$ and $T$ let
\[
\mathcal{P}(G,T)=\left\{P\in\mathcal{P}(G)\mid G/P\cong T\right\}
\]
be the collection of partitions of $G$ with quotient isomorphic to $T$. For
$P\in\mathcal{P}(G)$ let
\[
\left\Vert P\right\Vert =
{\displaystyle\prod\limits_{p\in P}}
(\vert p\vert-1)!
\]
where $p\in P$ runs through the $P$-classes.

Let $\mathrm{Aut}(G)$ denote the automorphism group of $G$. Let $\mathcal{G}
(k)$ denote the set of graphs without isolated vertices, where the number of
vertices minus the number of components equals $k$ and let
\[
\mathcal{G}(\leq k)=\cup_{j\leq k}\mathcal{G}(j)\text{.}
\]
Note that $\mathcal{G}(\leq k)$ is a finite set.\bigskip

\noindent\textbf{Base parameters.} Now we introduce a sequence of parameters
that will connect moments with chromatic coefficients. For a simple graph $T$
and $k>0$ let

\[
c_{k}(T)=
{\displaystyle\sum\limits_{G\in\mathcal{G}(k)}}
\frac{(-1)^{\left\vert E(G)\right\vert +\left\vert V(G)\right\vert +\left\vert
V(T)\right\vert +k}}{\left\vert \mathrm{Aut}(G)\right\vert}
{\displaystyle\sum\limits_{P\in\mathcal{P}(G,T)}}
\left\Vert P\right\Vert\text{.}
\]

It turns out that these parameters allow us to express $e_{k}(H)$ in a nice way.

\begin{lemma}
\label{ceka}Let $H$ be a simple graph. Then we have
\[
e_{k}(H)=
{\displaystyle\sum\limits_{T\in\mathcal{G}(\leq k)}}
c_{k}(T)\hom(T,H)\text{.}
\]

\end{lemma}

\noindent\textbf{Proof.} We derive the lemma from two easy claims. Let
$\mathrm{inj}(G,H)$ denote the number of injective homomorphisms from $G$ to
$H$.

\noindent\textbf{Claim 1.} We have
\[
e_{k}(H)=
{\displaystyle\sum\limits_{G\in\mathcal{G}(k)}}
\frac{(-1)^{\left\vert E(G)\right\vert +k}}{\left\vert \mathrm{Aut}
(G)\right\vert}\mathrm{inj}(G,H)\text{.}
\]
To see this, we use the following identity, that is sometimes used as a definition.

\[
\mathrm{ch}_{H}(z)=\sum_{\substack{G\subseteq H \\\text{spanning}
}}(-1)^{|E(G)|}z^{c(G)}
\]
where $c(G)$ is the number of connected components in the spanning subgraph
$G$. It is enough to prove this for positive integer values of $z$. In this
case, there are exactly $z^{c(G)}$ colorings that violate the legal coloring
constraint for all edges of $G$, and the equation follows from the
inclusion-exclusion principle.

The value of $e_{k}(H)$ is $(-1)^{k}$ times the coefficient of $z^{n-k}$,
which contains the terms where $c(G)=n-k$, or equivalently, where the graph
$G$, when erasing its isolated vertices, is in $\mathcal{G}_{k}$. A graph $G$
is counted as many times as it appears in $H$ as a spanning subgraph, which
equals $\mathrm{inj}(G,H)/\left\vert \mathrm{Aut}(G)\right\vert $. Claim 1 is proved.

\noindent\textbf{Claim 2.} Let $G\in\mathcal{G}(k)$ and let $H$ be a simple
graph. Then we have
\[
\mathrm{inj}(G,H)=
{\displaystyle\sum\limits_{T\in\mathcal{G}(\leq k)}}
\left(\left(-1\right)^{\left\vert V(G)\right\vert +\left\vert
V(T)\right\vert}
{\displaystyle\sum\limits_{P\in\mathcal{P}(G,T)}}
\left\Vert P\right\Vert \right)\hom(T,H)\text{.}
\]
To see this, let us consider the partially ordered set $\mathcal{P}(G)$ with
respect to refinement. For $P,P^{\prime}\in\mathcal{P}(G)$ with
$P^{\prime}\leq P$ (i.e. $P^{\prime}$ refines $P$), let $p_{1},\ldots,p_{r}$
be a list of the $P$-classes and let $a_{i}$ be the number of $P^{\prime}$-classes
contained in $p_{i}$ ($1\leq i\leq r$). Let $s=\displaystyle\sum_{i=1}^{r}a_{i}$
be the number of classes in $P^{\prime}$. Then the Mobius function is
\[
\mu(P^{\prime},P)=(-1)^{r+s}\prod_{i=1}^{r}(a_{i}-1)!
\]
(see e.g.\ \cite{rota}).
In particular, for the discrete partition $P_{0}$ we get
\[
\mu(P_{0},P)=(-1)^{|V(G)|+|V(G/P)|}\left\Vert P\right\Vert \text{.}
\]
On the other hand, we have
\[
\mathrm{\hom}(G/P^{\prime},H)=\sum_{P^{\prime}\leq P\in\mathcal{P}(G)}
\mathrm{inj}(G/P,H)\text{.}
\]
Now the Mobius inversion formula yields
\[
\mathrm{inj}(G,H)=\sum_{P\in\mathcal{P}(G)}(-1)^{|V(G)|+|V(G/P)|}\left\Vert
P\right\Vert \hom(G/P,H)
\]
which, when collecting terms by $T=G/P\in\mathcal{G}(\leq k)$ gives the
formula in Claim 2.

The lemma follows from substituting the formula in Claim 2 into the formula in
Claim 1 and collecting terms. $\square$\bigskip

Now we show that the base parameters of a disconnected graph can be expressed
as a convolution of the base parameters of its connected components,
normalized by a constant computed from the multiplicities:

\begin{lemma}
\label{convolution}Let $T$ be the disjoint union of the connected graphs
$T_{1},T_{2},\ldots,T_{l}$. Let $S = \{j \mid \nexists i<j: T_i \cong T_j\}$
contain the indices of nonisomorphic $T_j$'s and
$m_j = \left|\{i \mid T_i \cong T_j\}\right|$ denote the multiplicity of $T_j$.
Define ${\sigma=\displaystyle\prod_{j\in S} m_j!}$.
Then for all $k>0$ we have
\[
c_{k}(T)=\frac{1}{\sigma}
{\displaystyle\sum\limits_{\substack{(x_{1},\ldots,x_{l}) \\x_{1}+\cdots
+x_{l}=k}}}
{\displaystyle\prod\limits_{j=1}^{l}}
c_{x_{j}}(T_{j})\text{.}
\]
\end{lemma}

\noindent\textbf{Proof.} Recall that
\[
c_k(T)=\sum_{G\in\mathcal{G}(k)}
\frac{(-1)^{|E(G)|+|V(G)|+|V(T)|+k}}{|\mathrm{Aut}(G)|}
\sum_{P\in\mathcal{P}(G,T)} \left\Vert P\right\Vert \text{.}
\]

For a fixed $G$ and $P$, the connected components of $G/P$ can be identified
with $T_1,T_2,\ldots,T_l$ in $\sigma$ possible ways. Each of these matchings
gives a subdivision $G=G_1\cup G_2\cup\ldots\cup G_l$ by applying the inverse
image of the quotient map to the $T_i$'s. The restrictions $P_i=P|_{G_i}$ of the
partition $P$ satisfy $P_i\in\mathcal{P}(G_i,T_i)$ and
\[
\prod_{j=1}^l \left\Vert P_j\right\Vert = \left\Vert P\right\Vert \text{.}
\]
Therefore
\[
c_k(T) =
\sum_{G\in\mathcal{G}(k)}
\frac{(-1)^{|E(G)|+|V(G)|+|V(T)|+k}}{|\mathrm{Aut}(G)|}
\sum_{G=G_1\cup\ldots\cup G_l}
\sum_{\substack{P_j\in\mathcal{P}(G_j,T_j)\\1\le j\le l}}
\frac{1}{\sigma}
\prod_{j=1}^l \left\Vert P_j\right\Vert \text{.}
\]

If we already know the isomorphism classes of $G_1,G_2,\ldots,G_l$, there are
still
\[
\frac{|\mathrm{Aut}(G)|}{\displaystyle\prod_{j=1}^l |\mathrm{Aut}(G_j)|}
\]
possibilities to arrange them as a subdivision of $G$. It follows that
$c_k(T)$ equals
\[
\sum_{\substack{(x_1,\ldots,x_l)\\x_1+\ldots+x_l=k}}
\sum_{\substack{G_j\in\mathcal{G}(x_j)\\1\le j\le l}} \hspace{-2pt}
\frac{|\mathrm{Aut}(G)|}{\displaystyle\prod_{j=1}^l |\mathrm{Aut}(G_j)|}
\cdot
\frac{(-1)^{|E(G)|+|V(G)|+|V(T)|+k}}{|\mathrm{Aut}(G)|} \hspace{-12pt}
\sum_{\substack{P_j\in\mathcal{P}(G_j,T_j)\\1\le j\le l}} \hspace{-2pt}
\frac{1}{\sigma} \hspace{-1pt}
\prod_{j=1}^l \hspace{-1pt}
\left\Vert P_j\right\Vert \text{.}
\]

By using
\[
|E(G)|+|V(G)|+|V(T)|+k = \sum_{j=1}^l \left(|E(G_j)|+|V(G_j)|+|V(T_j)|+x_j\right)
\]
and rearranging we get
\begin{gather*}
c_k(T) =
\frac{1}{\sigma} \hspace{-5pt}
\sum_{\substack{(x_1,\ldots,x_l)\\x_1+\ldots+x_l=k}}
\prod_{j=1}^l
\sum_{G_j\in\mathcal{G}(x_j)} \hspace{-5pt}
\frac{(-1)^{|E(G_j)|+|V(G_j)|+|V(T_j)|+x_j}}{|\mathrm{Aut}(G_j)|}
\hspace{-10pt}
\sum_{P_j\in\mathcal{P}(G_j,T_j)}
\left\Vert P_j\right\Vert
= \\ =
\frac{1}{\sigma} \hspace{-5pt}
\sum_{\substack{(x_1,\ldots,x_l)\\x_1+\ldots+x_l=k}}
\prod_{j=1}^l c_{x_j}(T_j) \text{. $\square$}
\end{gather*}

We can use the following variant of Lemma \ref{convolution} when we would like
to detach one connected component of $T$ at a time:

\begin{lemma}
\label{detach}
Let $T$ be the disjoint union of the connected graphs $T_1,T_2,\ldots,T_l$ where
$l\ge 2$. Let $S = \{j \mid \nexists i<j: T_i \cong T_j\}$ contain the indices
of nonisomorphic $T_j$'s. Then we have
\[
kc_k(T)-
\sum_{i=1}^{k-1}
\sum_{j\in S}
ic_i(T_j)c_{k-i}(T\setminus T_j)
=0 \text{.}
\]
\end{lemma}

\noindent\textbf{Proof.}
Let $m_j$ denote the multiplicity of $T_j$ and
$\sigma=\displaystyle\prod_{j\in S} m_j!$ as in \mbox{Lemma \ref{convolution}}.
Since isomorphic $T_j$'s have identical $c_i(T_j)$ and $c_{k-i}(T\setminus T_j)$,
it follows that
\[
\sum_{j\in S}
ic_i(T_j)c_{k-i}(T\setminus T_j) =
\sum_{t=1}^l
\frac{i}{m_t}
c_i(T_t)c_{k-i}(T\setminus T_t) \text{.}
\]
By using Lemma \ref{convolution} for $T$ and $\sigma$ and also for
$T\setminus T_t$ and $\frac{\sigma}{m_t}$ we obtain:
\begin{gather*}
kc_k(T)-\hspace{-1pt}
\sum_{i=1}^{k-1}
\sum_{j\in S}
ic_i(T_j)c_{k-i}(T\setminus T_j) =
kc_k(T)-\hspace{-1pt}
\sum_{i=1}^{k-1}
\sum_{t=1}^l
\frac{i}{m_t}
c_i(T_t)c_{k-i}(T\setminus T_t)
= \\ =
\frac{k}{\sigma}\hspace{-4pt}
\sum_{\substack{(x_1,\ldots,x_l)\\x_1+\ldots+x_l=k}}
\prod_{j=1}^l
c_{x_j}(T_j)-
\sum_{i=1}^{k-1}
\sum_{t=1}^l
\frac{i}{m_t}\cdot
\frac{m_t}{\sigma}\hspace{-7pt}
\sum_{\substack{(x_1,\ldots,x_l)\\x_1+\ldots+x_l=k\\x_t=i}}
\prod_{j=1}^l
c_{x_j}(T_j)
= \\ =
\frac{1}{\sigma}\hspace{-5pt}
\sum_{\substack{(x_1,\ldots,x_l)\\x_1+\ldots+x_l=k}}
\left(k-\sum_{t=1}^l x_t\right)
\prod_{j=1}^l
c_{x_j}(T_j)
= 0 \text{.}
\end{gather*}

The last equation follows from
$k-\displaystyle\sum_{t=1}^l x_t=0$. $\square$

Now we show that $p_{k}(H)$ can be expressed using the number of homomorphisms
from \emph{connected} graphs.

\begin{theorem}
\label{fotetel}Let $H$ be a simple graph on $n$ vertices and let $k>0$ be an
integer. Then
\[
p_{k}(H)=
{\displaystyle\sum\limits_{\substack{T\in\mathcal{G}(\leq k)\\T\text{ is
connected}}}}
(-1)^{k-1}kc_{k}(T)\hom(T,H)\text{.}
\]

\end{theorem}

\noindent\textbf{Proof.} The Newton identites tell us that

\[
p_{k}(H)=(-1)^{k-1}ke_{k}(H)+
{\displaystyle\sum\limits_{i=1}^{k-1}}
(-1)^{k-i-1}p_{i}(H)e_{k-i}(H)
\]
for all $k>0$. Using induction on $k$, we can assume that the result holds for
all $j<k$. Lemma \ref{ceka} gives us a formula for $e_{k}(H)$ in the
parameters $c_{k}(T)$, namely we have
\[
e_{k}(H)=
{\displaystyle\sum\limits_{T\in\mathcal{G}(\leq k)}}
c_{k}(T)\hom(T,H)\text{.}
\]
\ 

Using that $\hom$ is multiplicative (stated as Lemma \ref{szorzat}) we get
$p_{k}(H)$ as a fixed linear combination of the $\hom(T,H)$'s. Let $q_{k}(T)$
denote the formal coefficient of $\hom(T,H)$ in this sum. So, we have
\[
p_{k}(H)=
{\displaystyle\sum\limits_{T\in\mathcal{G}(\leq k)}}
q_{k}(T)\hom(T,H)\text{.}
\]

This leads to the following equality for all $T$:
\[
q_{k}(T)=(-1)^{k-1}kc_{k}(T)+
{\displaystyle\sum\limits_{i=1}^{k-1}}
(-1)^{k-i-1}
{\displaystyle\sum\limits_{\substack{U_{1}\in\mathcal{G}(\leq i)\\U_{2}
\in\mathcal{G}(\leq k-i)\\U_{1}\cup U_{2}=T}}}
q_{i}(U_{1})c_{k-i}(U_{2})
\]
where $T$ is isomorphic to the disjoint union of $U_{1}$ and $U_{2}$.

Let $T\in\mathcal{G}(\leq k)$. We claim that
\[
q_{k}(T)=(-1)^{k-1}kc_{k}(T)
\]
if $T$ is connected and $0$ otherwise. If $T$ is connected then it is
impossible to choose $U_{1}$ and $U_{2}$ in the second sum above, so the claim
holds. If $T$ is disconnected then as in Lemma \ref{detach}, let $T$ be
the disjoint union of the connected graphs $T_{1},T_{2},\ldots,T_{l}$
and let $S$ contain the indices of nonisomorphic $T_j$'s.
Using induction on $k$ we can assume that $q_{i}(U_{1})=0$ unless
$U_1$ is isomorphic to one of the $T_j$'s. This gives
\[
q_{k}(T)=
(-1)^{k-1}kc_{k}(T)+
\sum_{i=1}^{k-1}
(-1)^{k-i-1}
\sum_{j\in S}
q_i(T_j)c_{k-i}(T\setminus T_j) \text{.}
\]

We know from the induction hypothesis that
$q_i(T_j)=(-1)^{i-1}ic_i(T_j)$ and therefore we get
\[
q_{k}(T)=
(-1)^{k-1}
\left(
kc_{k}(T)-
\sum_{i=1}^{k-1}
\sum_{j\in S}
ic_i(T_j)c_{k-i}(T\setminus T_j)
\right)
\]
which is 0 according to Lemma \ref{detach}. $\square$
\bigskip

\newpage
\section{Convergence of chromatic measures}
\label{convsection}

In this section we prove Theorem \ref{moments}, Theorem \ref{holom} and
Proposition \ref{tube}. For the convenience of the reader, we state the
theorems again.\bigskip

\begingroup
\def\thetheorem{\ref{moments}}
\begin{theorem}
Let $(G_{n})$ be a Benjamini-Schramm convergent graph sequence
of absolute degree bound $d$, and $\widetilde{D}$ an open neighborhood of
the closed disc $\overline{D}$. Then for every holomorphic function
$f:\widetilde{D}\rightarrow\mathbb{C}$, the sequence
\[
{\displaystyle\int\limits_{D}}
f(z)d\mu_{G_{n}}(z)
\]
converges.
\end{theorem}
\addtocounter{theorem}{-1}
\endgroup

\noindent\textbf{Proof.} We have
\[
\int\limits_{D}z^{k}d\mu_{G}(z)=\frac{1}{\left\vert V(G)\right\vert}
{\displaystyle\sum\limits_{i=1}^{\left\vert V(G)\right\vert}}
\lambda_{i}^{k}(G)=\frac{p_{k}(G)}{|V(G)|}
\]
for $k\geq0$.

Since $f$ is holomorphic, it equals its Taylor series
\[
f(z)=\sum_{n=0}^{\infty}a_{n}z^{n}
\]
on an open neighborhood of $\overline{D}$. Let
\[
f_k(z)=\sum_{n=0}^{k}a_{n}z^{n}
\]
denote the partial sums. The $f_k$'s converge uniformly on $D$, so
we also know that
\[
F_k(G)=
{\displaystyle\int\limits_{D}}
f_k(z)d\mu_{G}(z)=\sum_{n=0}^{k}a_{n}
{\displaystyle\int\limits_{D}}
z^{n}d\mu_{G}(z)=\sum_{n=0}^{k}a_{n}\frac{p_{n}(G)}{|V(G)|}
\]
converges to
\[
F(G)={\displaystyle\int\limits_{D}} f(z)d\mu_{G}(z)
\]
uniformly on the set of graphs $G$ with $\mu_G$ supported on $D$.
By Theorem \ref{fotetel} we have
\[
p_{n}(G)=
{\displaystyle\sum\limits_{\substack{T\in\mathcal{G}(\leq n)\\T\text{ is
connected}}}}
(-1)^{n-1}nc_{n}(T)\hom(T,G)\text{.}
\]
By rearranging, this gives
\[
F_k(G)=\sum_{T}b_{k,T}\frac{\hom(T,G)}{|V(G)|}
\]
where $T$ runs through connected graphs on at most $k+1$ vertices.
Now let $G_{n}$ be a Benjamini-Schramm convergent sequence of
graphs. By Lemma \ref{convhom}, the sequences
\[
\frac{\hom(T,G_{n})}{|V(G_{n})|}
\]
converge for every connected $T$. (Note that for non-connected $T$ this is in
general false). This implies that $F_k(G_n)$ is convergent for every $k$.
Since we already know that $F_k(G_n)$ uniformly converges to $F(G_n)$ for
every $n$, we obtain that $F(G_n)$ is also convergent.
It also follows that
\[
F(G_{n},u)=
{\displaystyle\int\limits_{D}}
f(z+u)d\mu_{G_{n}}(z)
\]
uniformly converges to a holomorphic function in a neighborhood of $0$.
$\square$\bigskip

We are ready to prove Theorem \ref{holom}.

\begingroup
\def\thetheorem{\ref{holom}}
\begin{theorem}
Let $(G_{n})$ be a Benjamini-Schramm convergent graph sequence of
absolute degree bound $d$ with $\left\vert V(G_{n})\right\vert \rightarrow\infty
$. Then $\mathrm{t}_{G_{n}}(z)$ converges to a real analytic function on
$\mathbb{C}\setminus\overline{D}$.
\end{theorem}
\addtocounter{theorem}{-1}
\endgroup

\noindent\textbf{Proof.} The principal branch of the
complex logarithm function only takes values with an imaginary part in
$(-\pi,\pi]$. Therefore $\Im\mathrm{t}_{G_{n}}(z)$ is always in the interval
$\left(\frac{-\pi}{\left\vert V(G_{n})\right\vert},\frac{\pi}{\left\vert
V(G_{n})\right\vert}\right]$ and $\left\vert V(G_{n})\right\vert
\rightarrow\infty$ implies $\Im\mathrm{t}_{G_{n}}(z)\rightarrow0$.

To prove convergence for the real part $\Re\mathrm{t}_{G_{n}}(z)$, consider a
fixed $z_{0}\in\mathbb{C}\setminus\overline{D}$. Since the disc $B(z_{0},Cd)$
is bounded away from 0, there exists a branch $\ln^{\ast}$ of the complex
logarithm function whose branch cut is a half-line emanating from $0$ that is
disjoint from the disc. It follows that $f(z)=\ln^{\ast}(z_{0}-z)$ is holomorphic
on an open neighborhood of $\overline{D}$. According to Theorem \ref{moments},
\[
\displaystyle\int\limits_{D}\ln^{\ast}(z_{0}-z)d\mu_{G_{n}}(z)
\]
converges uniformly in a neighborhood of $z_{0}$, which implies that
\begin{gather*}
\Re\mathrm{t}_{G_{n}}(z_{0})=\frac{\Re\ln \mathrm{ch}_{G_{n}}(z_{0})}{\left\vert
V(G_{n})\right\vert}=\frac{\displaystyle\sum_{\text{$\lambda$ root}}\Re
\ln(z_{0}-\lambda)}{\left\vert V(G_{n})\right\vert}=\displaystyle\int
\limits_{D}\Re\ln(z_{0}-z)d\mu_{G_{n}}(z)=\\
=\displaystyle\int\limits_{D}\Re\ln^{\ast}(z_{0}-z)d\mu_{G_{n}}(z)=\Re
\int\limits_{D}\ln^{\ast}(z_{0}-z)d\mu_{G_{n}}(z)
\end{gather*}
is locally uniformly convergent as a function of $z_{0}$. Since $\Re\ln
(z_{0}-\lambda)$ is a harmonic function for all chromatic roots $\lambda$, so
is $\Re\mathrm{t}_{G_{n}}(z_{0})$, and the harmonicity of $\lim\mathrm{t}
_{G_{n}}(z_{0})=\lim\Re\mathrm{t}_{G_{n}}(z_{0})$ follows from local uniform
convergence. The observation that all harmonic functions are real analytic
concludes the proof. $\square$\bigskip

Now we prove Proposition \ref{tube}. Note that already Salas and Sokal
\cite{salas} showed that the pointwise limit of supports of $\mu_{T_{n}}$ is
part of a particular algebraic curve. For convenience, we include some details
on that, also adding a picture on the supporting set, but we do not introduce
the transfer matrix method here. See \cite{salas} for a description of the
transfer matrix method.

\begingroup
\def\thetheorem{\ref{tube}}
\begin{proposition}
The chromatic measures $\mu_{T_{n}}$ weakly converge.
\end{proposition}
\addtocounter{theorem}{-1}
\endgroup

\noindent\textbf{Proof.} We defined $T_{n}$ as the
cartesian product of $C_{4}$ and $P_{n}$. By the transfer matrix method we
obtain
\[
\mathrm{ch}_{T_{n}}(z)=v_{1}M^{n-1}\underline{\mathbf{{1}}}^{\mathsf{T}}
\]
with
\begin{gather*}
v_{1}=\left(
\begin{array}
[c]{ccc}
z^{4}{-}6z^{3}{+}11z^{2}{-}6z & 2z^{3}{-}6z^{2}{+}4z & z^{2}{-}z
\end{array}
\right)\\[5pt]
M=\left(
\begin{array}
[c]{ccc}
z^{4}{-}10z^{3}{+}41z^{2}{-}84z{+}73 & 2z^{3}{-}14z^{2}{+}38z{-}40 & z^{2}{-}5z{+}8\\
z^{4}{-}10z^{3}{+}40z^{2}{-}77z{+}60 & 2z^{3}{-}13z^{2}{+}32z{-}29 & z^{2}{-}4z{+}5\\
z^{4}{-}10z^{3}{+}39z^{2}{-}70z{+}48 & 2z^{3}{-}12z^{2}{+}26z{-}20 & z^{2}{-}3z{+}3
\end{array}
\right)\\[5pt]
\underline{\mathbf{{1}}}^{\mathsf{T}}=\left(
\begin{array}
[c]{c}
1\\
1\\
1
\end{array}
\right)\raisebox{-12pt}{.}
\end{gather*}

Using the eigenvectors of $M$ as our new basis we can diagonalize $M$ and
rewrite the above expression as
\[
\mathrm{ch}_{T_{n}}(z)=u_{1}D^{n-1}u_{2}
\]
where
\begin{gather*}
u_{1}=\left(
\begin{array}
[c]{c}
\frac{z^{7}{-}10z^{6}{+}44z^{5}{-}105z^{4}{+}143z^{3}{-}109z^{2}{+}36z{+}z^{3}r{-}2z^{2}
r{+}zr}{2z^{3}{-}12z^{2}{+}28z{-}24}\\
\frac{z^{7}{-}10z^{6}{+}44z^{5}{-}105z^{4}{+}143z^{3}{-}109z^{2}{+}36z{-}z^{3}r{+}2z^{2}
r{-}zr}{2z^{3}{-}12z^{2}{+}28z{-}24}\\
0
\end{array}
\right)^{\mathsf{T}}\\
D=\left(
\begin{array}
[c]{ccc}
\frac{z^{4}{-}8z^{3}{+}29z^{2}{-}55z{+}46{+}r}{2} & 0 & 0\\
0 & \frac{z^{4}{-}8z^{3}{+}29z^{2}{-}55z{+}46{-}r}{2} & 0\\
0 & 0 & 1
\end{array}
\right) \\
u_{2}=\left(
\begin{array}
[c]{c}
\frac{z^{4}{-}8z^{3}{+}27z^{2}{-}47z{+}36{+}r}{2r}\\
\frac{{-}z^{4}{+}8z^{3}{-}27z^{2}{+}47z{-}36{+}r}{2r}\\
0
\end{array}
\right)
\end{gather*}
and
\[
r=\sqrt{z^{8}{-}16z^{7}{+}118z^{6}{-}526z^{5}{+}1569z^{4}{-}3250z^{3}
{+}4617z^{2}{-}4136z{+}1776}\text{.}
\]

The matrix $D^{n-1}$ is straightforward to calculate, so we get the following closed
formula for the chromatic polynomial:
\[
\mathrm{ch}_{T_{n}}(z)=a_1\lambda_1^{n-1}+a_2\lambda_2^{n-1}
\]
where
\[
a_i=
\frac{z(z{-}1)(z^{4}{-}8z^{3}{+}27z^{2}{-}47z{+}36{+}r_{i})(z^{5}{-}9z^{4}{+}35z^{3}{-}70z^{2}{+}73z{-}36{+}zr_{i}{-}r_{i})}
{4r_{i}(z{-}2)(z^{2}{-}4z{+}6)}
\]
and
\[
\lambda_i=\frac{z^{4}{-}8z^{3}{+}29z^{2}{-}55z{+}46{+}r_{i}}{2}
\]
with $r_{1,2}=\pm r$.

We are interested in the complex roots of this expression if $n$ is very
large. We don't need to specify them exactly, but we'll prove a necessary
condition. If the eigenvalues $\lambda_{i}$
differ in their absolute value for some $z$, there will be
an arbitrarily large multiplicative gap between $a_{1}\lambda_{1}^{n-1}$ and
$a_{2}\lambda_{2}^{n-1}$ for any values of $a_{i}$ unless both $a_{1}
\lambda_{1}=0$ and $a_{2}\lambda_{2}=0$ holds.

It follows that all roots must have $|\lambda_{1}|=|\lambda_{2}|$ with the
possible exception of a finite set consisting of the roots and singularities
of $a_{1}$, $a_{2}$, $\lambda_{1}$ and $\lambda_{2}$, or equivalently, the
roots of
\begin{gather*}
z(z{-}1)(z{-}2)(2z{-}5)(z^{2}{-}3z{+}1)(z^{2}{-}4z{+}6)\cdot\\
\cdot(z^{6}{-}12z^{5}{+}61z^{4}{-}169z^{3}{+}269z^{2}{-}231z{+}85)\cdot\\
\cdot(z^{8}{-}16z^{7}{+}118z^{6}{-}526z^{5}{+}1569z^{4}{-}3250z^{3}{+}4617z^{2}
{-}4136z{+}1776)\text{.}
\end{gather*}

Let's ignore this set $\mathcal{S}$ of special roots for now and concentrate
on the general case of $|\lambda_{1}|=|\lambda_{2}|$:
\begin{gather*}
\lambda_{1}\overline{\lambda_{1}}=\lambda_{2}\overline{\lambda_{2}}\\
\frac{z^{4}{-}8z^{3}{+}29z^{2}{-}55z{+}46{+}r}{2}\cdot\frac{\overline{z}^{4}
{-}8\overline{z}^{3}{+}29\overline{z}^{2}{-}55\overline{z}{+}46{+}\overline{r}}{2}=\\
=\frac{z^{4}{-}8z^{3}{+}29z^{2}{-}55z{+}46{-}r}{2}\cdot\frac{\overline{z}^{4}
{-}8\overline{z}^{3}{+}29\overline{z}^{2}{-}55\overline{z}{+}46{-}\overline{r}}{2}\\
(z^{4}{-}8z^{3}{+}29z^{2}{-}55z{+}46)\overline{r}\hspace{1pt}+(\overline{z}^{4}{-}8\overline{z}
^{3}{+}29\overline{z}^{2}{-}55\overline{z}{+}46)r=0
\end{gather*}

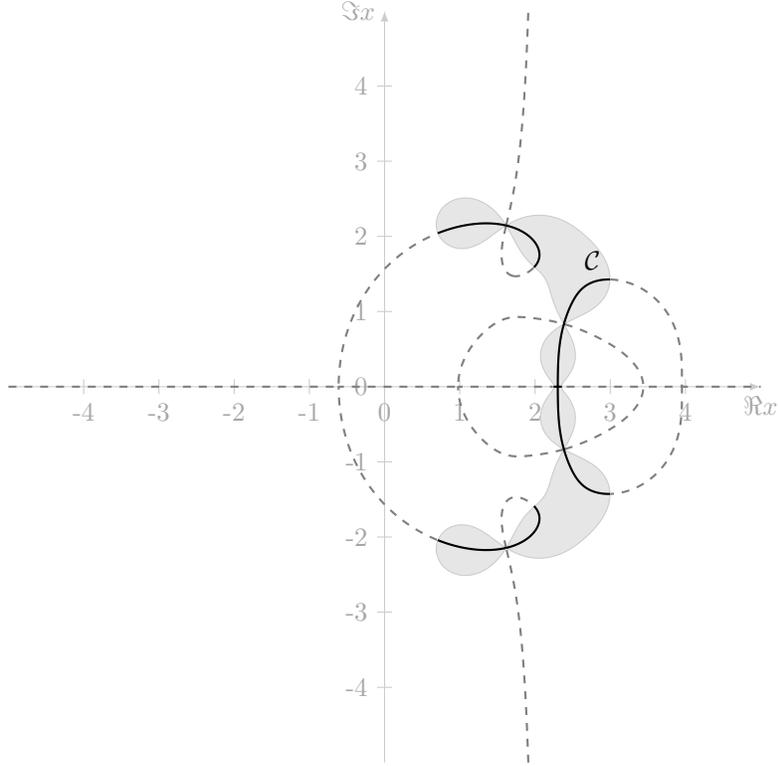
\begin{figure}[h]
\begin{center}
\tikzsetnextfilename{tube}
\begin{tikzpicture}
\draw[thin,color=black!20,fill=black!10]
plot file {plots/tubec00.table}
plot file {plots/tubec01.table}
plot file {plots/tubec02.table};
\draw[thin,color=black!20,-latex]
(-5,0) -- (5,0)
node[below,color=black!35] {$\Re x$};
\draw[thin,color=black!20,-latex]
(0,-5) -- (0,5)
node[left,color=black!35] {$\Im x$};
\foreach \i in {-4,...,4}
\draw[thin,color=black!20]
(\i,0.1) -- (\i,-0.1)
node[below,color=black!35] {\i};
\foreach \i in {-4,...,4}
\draw[thin,color=black!20]
(0.1,\i) -- (-0.1,\i)
node[left,color=black!35] {\i};
\draw[thick,style=dashed,color=black!50]
plot file {plots/tubee00.table}
plot file {plots/tubee01.table}
plot file {plots/tubee02.table}
plot file {plots/tubee03.table}
plot file {plots/tubee04.table}
plot file {plots/tubee05.table}
plot file {plots/tubee06.table}
plot file {plots/tubee07.table}
plot file {plots/tubee08.table}
plot file {plots/tubee09.table}
plot file {plots/tubee10.table};
\draw[thick,color=black]
plot file {plots/tuber00.table}
plot file {plots/tuber01.table}
node[above left] {$\mathcal{C}$}
plot file {plots/tuber02.table}
plot file {plots/tuber03.table}
plot file {plots/tuber04.table}
plot file {plots/tuber05.table}
plot file {plots/tuber06.table}
plot file {plots/tuber07.table}
plot file {plots/tuber08.table};
\end{tikzpicture}
\end{center}
\caption{Possible limit points of chromatic roots of $T_n=C_4\times P_n$ as $n\to\infty$}
\label{tubepicture}
\end{figure}

Our last expression means that $(z^{4}{-}8z^{3}{+}29z^{2}{-}55z{+}46)\overline{r}$ is
purely imaginary, which is equivalent to its square being a nonpositive real.
When calculated, this gives a degree 14 algebraic curve clipped by a degree 16
algebraic curve, shown as the set $\mathcal{C}$ on the Figure \ref{tubepicture}.

It follows that the curve $\mathcal{C}$ is compact, has an empty interior and
its complement is connected. Hence the same holds for $\mathcal{C}^{\prime
}=\mathcal{C\cup S}$.

Now Mergelyan's theorem \cite{merg} says that every continuous function on
$\mathcal{C}^{\prime}$ can be uniformly approximated by polynomials. This
implies that if two probability measures $\mu_{1}$ and $\mu_{2}$ are supported
on $\mathcal{C}^{\prime}$ and the holomorphic moments satisfy
\[
{\displaystyle\int\limits_{\mathcal{C}^{\prime}}}
z^{k}d\mu_{1}(z)=
{\displaystyle\int\limits_{\mathcal{C}^{\prime}}}
z^{k}d\mu_{2}(z)\text{ \ (}k\geq1\text{)}
\]
then
\[
{\displaystyle\int\limits_{\mathcal{C}^{\prime}}}
f(z)d\mu_{1}(z)=
{\displaystyle\int\limits_{\mathcal{C}^{\prime}}}
f(z)d\mu_{2}(z)
\]
for all continuous functions $f:\mathcal{C}^{\prime}\rightarrow\mathbb{R}$.
Hence, we have $\mu_{1}=\mu_{2}$. Since any subsequential weak limit of
$\mu_{T_{n}}$ is supported on $\mathcal{C}^{\prime}$, we get that $\mu_{T_{n}
}$ is weakly convergent. $\square$\bigskip

\noindent\textbf{Remark.} As we saw in the introduction, weak convergence does
not hold in general. The phenomenon where an associated measure blows up by a
small change of the graph but keeps its holomorphic moments unchanged also
occurs in the spectral theory of \emph{directed} graphs. Namely, the weak
limit of the eigenvalue distributions of the directed path of length $n$ is
the Dirac measure at $0$, while for the directed $n$-cycle the limit is the
uniform measure on the unit circle centered at $0$. In both the chromatic and
the spectral case, the reason is that the change only affects the coefficients
of small index in the corresponding polynomial, and the $k$-th moment only
depends on the $k$ highest index coefficients. It would be interesting to
study this blow-up phenomenon using just abstract polynomials.

\section{Graphs of large girth}

In this section we study graphs with large girth and prove Theorem
\ref{girth}.\bigskip

\begin{lemma}
Suppose that the finite graphs $G$ and $H$ both have girth at least $g$ and
$|E(G)|=|E(H)|$. Then $e_{i}(G)=e_{i}(H)$ holds for $i=0,1,\ldots,g-2$.
\end{lemma}

\noindent\textbf{Proof.} We use induction on $|E(G)|=|E(H)|$. If the number of
edges is zero, the claim is trivial, as is when $g\leq3$. Otherwise pick $e\in
E(G)$ and $f\in E(H)$ arbitrarily and use the deletion-contraction argument:
\begin{gather*}
e_{i}(G)=e_{i}(G\setminus e)+e_{i-1}(G/e)\\
e_{i}(H)=e_{i}(H\setminus f)+e_{i-1}(H/f)
\end{gather*}

The claim follows from the observation that $G\setminus e$ and $H\setminus f$
have girth $\geq g$ while $G/e$ and $H/f$ have girth $\geq g-1$. $\square
$\bigskip

\begin{lemma}
\label{girthmoments}Let $G$ be a finite graph with girth $\geq g$. Then
$p_{i}(G)=|E(G)|$ for $i=1,2,\ldots,g-2$.
\end{lemma}

\noindent\textbf{Proof.} Let $H$ be an arbitrary tree on $|E(G)|+1$ vertices
and use the previous lemma. Since the chromatic polynomial of $H$ is
$q(q-1)^{|E(G)|}$, we have $e_{i}(H)=\binom{|E(G)|}{i}$ for $i\leq|E(G)|$,
which translates into $e_{i}(G)=\binom{|E(G)|}{i}$ for $i\leq g-2$.
Substituting the $e_{i}$'s into Newton's identities completes the proof.
$\square$\bigskip

\begingroup
\def\thetheorem{\ref{girth}}
\begin{theorem}
Let $G$ be a finite graph of girth $g$ and maximal degree $d$.
Then for all $q>Cd$ we have
\[
\left\vert \frac{\ln \mathrm{ch}_{G}(q)}{\left\vert V(G)\right\vert}-\left(\ln
q+\frac{\left\vert E(G)\right\vert}{\left\vert V(G)\right\vert}\ln
(1-\frac{1}{q})\right)\right\vert \leq2\frac{(Cd/q)^{g-1}}{1-Cd/q}\text{.}
\]
\end{theorem}
\addtocounter{theorem}{-1}
\endgroup

\noindent\textbf{Proof.} The normalized log of the
chromatic polynomial can be expanded as
\begin{gather*}
\frac{\ln \mathrm{ch}_{G}(q)}{|V(G)|}=\int\limits_{D}\ln(q-z)d\mu_{G}(z)=\ln
q+\int\limits_{D}\ln\left(1-\frac{z}{q}\right)d\mu_{G}(z)=\\
=\ln q-\sum_{n=1}^{\infty}\frac{1}{nq^{n}}\int\limits_{D}z^{n}d\mu_{G}(z)
\end{gather*}
where the Sokal bound $|z|\leq Cd$ gives the constraint
\[
\left\vert \int\limits_{D}z^{n}d\mu_{G}(z)\right\vert \leq(Cd)^{n}
\]
for the holomorphic moments, and our last lemma implies
\[
\int\limits_{D}z^{n}d\mu_{G}(z)=\frac{p_{n}(G)}{|V(G)|}=\frac{|E(G)|}{|V(G)|}
\]
for $n\leq g-2$. We also know that any real number $x\in\lbrack0,1)$
satisfies
\[
\sum_{n=g-1}^{\infty}\frac{x^{n}}{n}=\int\limits_{0}^{x}\sum_{n=g-2}^{\infty
}t^{n}dt=\int\limits_{0}^{x}\frac{t^{g-2}}{1-t}dt\leq x\cdot\frac{x^{g-2}
}{1-x}=\frac{x^{g-1}}{1-x}\;\text{.}
\]

Now we have

\begin{gather*}
\left\vert \frac{\ln \mathrm{ch}_{G}(q)}{|V(G)|}-\left(\ln q+\frac{|E(G)|}{|V(G)|}
\ln(1-\frac{1}{q})\right)\right\vert =\\
=\left\vert \left(\ln q-\sum_{n=1}^{\infty}\frac{1}{nq^{n}}\int
\limits_{D}z^{n}d\mu_{G}(Z)\right)-\left(\ln q-\sum_{n=1}^{\infty}\frac
{1}{nq^{n}}\cdot\frac{|E(G)|}{|V(G)|}\right)\right\vert \leq\\
\leq\left\vert \sum_{n=g-1}^{\infty}\frac{1}{nq^{n}}\left(\int
\limits_{D}z^{n}d\mu_{G}(z)-\frac{|E(G)|}{|V(G)|}\right)\right\vert \leq
\sum_{n=g-1}^{\infty}\frac{1}{nq^{n}}\left(\left\vert \int\limits_{D}
z^{n}d\mu_{G}(z)\right\vert +\frac{|E(G)|}{|V(G)|}\right)\leq\\
\leq\sum_{n=g-1}^{\infty}\frac{(Cd)^{n}+|E(G)|/|V(G)|}{nq^{n}}\leq\sum
_{n=g-1}^{\infty}\frac{2(Cd)^{n}}{nq^{n}}\leq2\frac{(Cd/q)^{g-1}}
{1-Cd/q}\;\text{.}
\end{gather*}

The theorem holds. $\square$

\newpage
\section{Appendix}

In the appendix we publish some data that may be useful for further analysis.

For abbreviation, we use the following terminology:
\[
\hom\left(
{\displaystyle\sum\limits_{i=1}^{n}}
\alpha_{i}G_{i},H\right)=
{\displaystyle\sum_{i=1}^{n}}
\alpha_{i}\hom(G_{i},H)
\]
where the $G_{i}$ and $H$ are finite graphs.

One can express the first $4$ chromatic coefficients as a linear combination
of homomorphisms as follows:\ \bigskip

$e_{0}(G)=\hom\hspace{-2pt}\left(\raisebox{0pt}[3pt][3pt]{}\right.
\hspace{-2pt}\ensuremath{\tikzsetfigurename{e0-}
\raisebox{0.5ex}{\tikzsetnextfilename{e0-0}\begin{tikzpicture}[scale=0.150000000000000,baseline=(current bounding box.center),n/.style={fill,shape=circle,inner sep=0pt,minimum size=2.50000000000000pt},>=latex,line join=bevel,]
\node (0) at (7bp,7.5bp) [n] {};
\end{tikzpicture}}
}
,G\hspace{-2pt}\left.
\raisebox{0pt}[3pt][3pt]{}\right)\hspace{-2pt}$

\bigskip

$e_{1}(G)=\hom\hspace{-2pt}\left(\raisebox{0pt}[4pt][4pt]{}\right.
\hspace{-3pt}\ensuremath{\tikzsetfigurename{e1-}
\frac{1}{2}\raisebox{0.5ex}{\tikzsetnextfilename{e1-0}\begin{tikzpicture}[scale=0.150000000000000,baseline=(current bounding box.center),n/.style={fill,shape=circle,inner sep=0pt,minimum size=2.50000000000000pt},>=latex,line join=bevel,]
\node (1) at (78.21bp,7.5bp) [n] {};
  \node (0) at (7bp,18.139bp) [n] {};
  \draw [] (0) ..controls (26.525bp,15.222bp) and (58.714bp,10.413bp)  .. (1);
\end{tikzpicture}}
},G\hspace{-2pt}\left.
\raisebox{0pt}[4pt][4pt]{}\right)\hspace{-2pt}$

\bigskip

$e_{2}(G)=\hom\hspace{-2pt}\left(\raisebox{0pt}[4pt][4pt]{}\right.
\hspace{-3pt}\ensuremath{\tikzsetfigurename{e2-}
-\frac{1}{4}\raisebox{0.5ex}{\tikzsetnextfilename{e2-0}\begin{tikzpicture}[scale=0.150000000000000,baseline=(current bounding box.center),n/.style={fill,shape=circle,inner sep=0pt,minimum size=2.50000000000000pt},>=latex,line join=bevel,]
\node (1) at (78.21bp,7.5bp) [n] {};
  \node (0) at (7bp,18.139bp) [n] {};
  \draw [] (0) ..controls (26.525bp,15.222bp) and (58.714bp,10.413bp)  .. (1);
\end{tikzpicture}} - \frac{1}{6}\raisebox{0.5ex}{\tikzsetnextfilename{e2-1}\begin{tikzpicture}[scale=0.150000000000000,baseline=(current bounding box.center),n/.style={fill,shape=circle,inner sep=0pt,minimum size=2.50000000000000pt},>=latex,line join=bevel,]
\node (1) at (79.206bp,69.908bp) [n] {};
  \node (0) at (7bp,69.434bp) [n] {};
  \node (2) at (43.592bp,7.5bp) [n] {};
  \draw [] (0) ..controls (26.798bp,69.564bp) and (59.438bp,69.778bp)  .. (1);
  \draw [] (0) ..controls (17.922bp,50.948bp) and (32.528bp,26.227bp)  .. (2);
  \draw [] (1) ..controls (68.576bp,51.28bp) and (54.361bp,26.37bp)  .. (2);
\end{tikzpicture}} + \frac{1}{8}\raisebox{0.5ex}{\tikzsetnextfilename{e2-2}\begin{tikzpicture}[scale=0.150000000000000,baseline=(current bounding box.center),n/.style={fill,shape=circle,inner sep=0pt,minimum size=2.50000000000000pt},>=latex,line join=bevel,]
\node (1) at (78.21bp,7.5bp) [n] {};
  \node (0) at (7bp,18.139bp) [n] {};
  \node (3) at (134.21bp,35.5bp) [n] {};
  \node (2) at (63bp,46.139bp) [n] {};
  \draw [] (0) ..controls (26.525bp,15.222bp) and (58.714bp,10.413bp)  .. (1);
  \draw [] (2) ..controls (82.525bp,43.222bp) and (114.71bp,38.413bp)  .. (3);
\end{tikzpicture}}
},G\hspace{-2pt}\left.
\raisebox{0pt}[4pt][4pt]{}\right)\hspace{-2pt}$

\bigskip

$e_{3}(G)=\hom\hspace{-2pt}\left(\raisebox{0pt}[6pt][6pt]{}\right.
\hspace{-4pt}\ensuremath{\tikzsetfigurename{e3-}
\frac{1}{24}\raisebox{0.5ex}{\tikzsetnextfilename{e3-0}\begin{tikzpicture}[scale=0.150000000000000,baseline=(current bounding box.center),n/.style={fill,shape=circle,inner sep=0pt,minimum size=2.50000000000000pt},>=latex,line join=bevel,]
\node (1) at (78.21bp,7.5bp) [n] {};
  \node (0) at (7bp,18.139bp) [n] {};
  \draw [] (0) ..controls (26.525bp,15.222bp) and (58.714bp,10.413bp)  .. (1);
\end{tikzpicture}} + \frac{1}{4}\raisebox{0.5ex}{\tikzsetnextfilename{e3-1}\begin{tikzpicture}[scale=0.150000000000000,baseline=(current bounding box.center),n/.style={fill,shape=circle,inner sep=0pt,minimum size=2.50000000000000pt},>=latex,line join=bevel,]
\node (1) at (149.73bp,20.115bp) [n] {};
  \node (0) at (7bp,17.213bp) [n] {};
  \node (2) at (78.597bp,7.5bp) [n] {};
  \draw [] (1) ..controls (130.23bp,16.656bp) and (98.072bp,10.954bp)  .. (2);
  \draw [] (0) ..controls (26.631bp,14.55bp) and (58.995bp,10.159bp)  .. (2);
\end{tikzpicture}} + \frac{1}{12}\raisebox{0.5ex}{\tikzsetnextfilename{e3-2}\begin{tikzpicture}[scale=0.150000000000000,baseline=(current bounding box.center),n/.style={fill,shape=circle,inner sep=0pt,minimum size=2.50000000000000pt},>=latex,line join=bevel,]
\node (1) at (79.206bp,69.908bp) [n] {};
  \node (0) at (7bp,69.434bp) [n] {};
  \node (2) at (43.592bp,7.5bp) [n] {};
  \draw [] (0) ..controls (26.798bp,69.564bp) and (59.438bp,69.778bp)  .. (1);
  \draw [] (0) ..controls (17.922bp,50.948bp) and (32.528bp,26.227bp)  .. (2);
  \draw [] (1) ..controls (68.576bp,51.28bp) and (54.361bp,26.37bp)  .. (2);
\end{tikzpicture}} - \frac{1}{8}\raisebox{0.5ex}{\tikzsetnextfilename{e3-3}\begin{tikzpicture}[scale=0.150000000000000,baseline=(current bounding box.center),n/.style={fill,shape=circle,inner sep=0pt,minimum size=2.50000000000000pt},>=latex,line join=bevel,]
\node (1) at (78.21bp,7.5bp) [n] {};
  \node (0) at (7bp,18.139bp) [n] {};
  \node (3) at (134.21bp,35.5bp) [n] {};
  \node (2) at (63bp,46.139bp) [n] {};
  \draw [] (0) ..controls (26.525bp,15.222bp) and (58.714bp,10.413bp)  .. (1);
  \draw [] (2) ..controls (82.525bp,43.222bp) and (114.71bp,38.413bp)  .. (3);
\end{tikzpicture}} - \frac{1}{8}\raisebox{0.5ex}{\tikzsetnextfilename{e3-4}\begin{tikzpicture}[scale=0.150000000000000,baseline=(current bounding box.center),n/.style={fill,shape=circle,inner sep=0pt,minimum size=2.50000000000000pt},>=latex,line join=bevel,]
\node (1) at (98.357bp,65.742bp) [n] {};
  \node (0) at (7bp,43.092bp) [n] {};
  \node (3) at (52.478bp,101.75bp) [n] {};
  \node (2) at (53.166bp,7.5bp) [n] {};
  \draw [] (0) ..controls (30.049bp,48.807bp) and (75.732bp,60.133bp)  .. (1);
  \draw [] (0) ..controls (21.834bp,31.656bp) and (38.099bp,19.116bp)  .. (2);
  \draw [] (1) ..controls (83.114bp,77.705bp) and (67.187bp,90.206bp)  .. (3);
  \draw [] (2) ..controls (52.989bp,31.668bp) and (52.653bp,77.711bp)  .. (3);
\end{tikzpicture}} + \frac{1}{4}\raisebox{0.5ex}{\tikzsetnextfilename{e3-5}\begin{tikzpicture}[scale=0.150000000000000,baseline=(current bounding box.center),n/.style={fill,shape=circle,inner sep=0pt,minimum size=2.50000000000000pt},>=latex,line join=bevel,]
\node (1) at (135.52bp,52.44bp) [n] {};
  \node (0) at (7bp,31.674bp) [n] {};
  \node (3) at (65.649bp,76.768bp) [n] {};
  \node (2) at (76.843bp,7.5bp) [n] {};
  \draw [] (0) ..controls (24.624bp,45.224bp) and (48.41bp,63.513bp)  .. (3);
  \draw [] (1) ..controls (117.89bp,38.936bp) and (94.09bp,20.71bp)  .. (2);
  \draw [] (0) ..controls (26.818bp,24.814bp) and (57.43bp,14.219bp)  .. (2);
  \draw [] (1) ..controls (115.69bp,59.343bp) and (85.07bp,70.006bp)  .. (3);
  \draw [] (2) ..controls (73.625bp,27.415bp) and (68.832bp,57.074bp)  .. (3);
\end{tikzpicture}} - \frac{1}{24}\raisebox{0.5ex}{\tikzsetnextfilename{e3-6}\begin{tikzpicture}[scale=0.150000000000000,baseline=(current bounding box.center),n/.style={fill,shape=circle,inner sep=0pt,minimum size=2.50000000000000pt},>=latex,line join=bevel,]
\node (1) at (92.936bp,56.852bp) [n] {};
  \node (0) at (7bp,44.132bp) [n] {};
  \node (3) at (43.555bp,93.513bp) [n] {};
  \node (2) at (56.291bp,7.5bp) [n] {};
  \draw [] (0) ..controls (19.067bp,60.432bp) and (31.554bp,77.301bp)  .. (3);
  \draw [] (2) ..controls (52.923bp,30.249bp) and (46.985bp,70.35bp)  .. (3);
  \draw [] (0) ..controls (22.504bp,32.609bp) and (40.707bp,19.082bp)  .. (2);
  \draw [] (1) ..controls (77.404bp,68.384bp) and (59.168bp,81.922bp)  .. (3);
  \draw [] (0) ..controls (29.036bp,47.394bp) and (71.018bp,53.608bp)  .. (1);
  \draw [] (1) ..controls (80.84bp,40.561bp) and (68.322bp,23.702bp)  .. (2);
\end{tikzpicture}} - \frac{1}{12}\raisebox{0.5ex}{\tikzsetnextfilename{e3-7}\begin{tikzpicture}[scale=0.150000000000000,baseline=(current bounding box.center),n/.style={fill,shape=circle,inner sep=0pt,minimum size=2.50000000000000pt},>=latex,line join=bevel,]
\node (1) at (79.206bp,69.908bp) [n] {};
  \node (0) at (7bp,69.434bp) [n] {};
  \node (3) at (88bp,36.139bp) [n] {};
  \node (2) at (43.592bp,7.5bp) [n] {};
  \node (4) at (159.21bp,25.5bp) [n] {};
  \draw [] (0) ..controls (26.798bp,69.564bp) and (59.438bp,69.778bp)  .. (1);
  \draw [] (0) ..controls (17.922bp,50.948bp) and (32.528bp,26.227bp)  .. (2);
  \draw [] (1) ..controls (68.576bp,51.28bp) and (54.361bp,26.37bp)  .. (2);
  \draw [] (3) ..controls (107.52bp,33.222bp) and (139.71bp,28.413bp)  .. (4);
\end{tikzpicture}} + \frac{1}{48}\raisebox{0.5ex}{\tikzsetnextfilename{e3-8}\begin{tikzpicture}[scale=0.150000000000000,baseline=(current bounding box.center),n/.style={fill,shape=circle,inner sep=0pt,minimum size=2.50000000000000pt},>=latex,line join=bevel,]
\node (1) at (78.21bp,7.5bp) [n] {};
  \node (0) at (7bp,18.139bp) [n] {};
  \node (3) at (134.21bp,35.5bp) [n] {};
  \node (2) at (63bp,46.139bp) [n] {};
  \node (5) at (99.21bp,63.5bp) [n] {};
  \node (4) at (28bp,74.139bp) [n] {};
  \draw [] (0) ..controls (26.525bp,15.222bp) and (58.714bp,10.413bp)  .. (1);
  \draw [] (4) ..controls (47.525bp,71.222bp) and (79.714bp,66.413bp)  .. (5);
  \draw [] (2) ..controls (82.525bp,43.222bp) and (114.71bp,38.413bp)  .. (3);
\end{tikzpicture}}
},G\hspace{-2pt}\left.
\raisebox{0pt}[6pt][6pt]{}\right)\hspace{-2pt}$

\bigskip

$e_{4}(G)=\hom\hspace{-2pt}\left(\raisebox{0pt}[8pt][8pt]{}\right.
\hspace{-5pt}\input{includes/e4.inc},G\hspace{-3pt}\left.
\raisebox{0pt}[8pt][8pt]{}\right)\hspace{-2pt}$

\bigskip

Also, one can express the first $5$ chromatic moments as a linear combination
of homomorphisms as follows. \bigskip

$p_{0}(G)=\hom\hspace{-2pt}\left(\raisebox{0pt}[3pt][3pt]{}\right.
\hspace{-2pt}\ensuremath{\tikzsetfigurename{p0-}
\raisebox{0.5ex}{\tikzsetnextfilename{p0-0}\begin{tikzpicture}[scale=0.150000000000000,baseline=(current bounding box.center),n/.style={fill,shape=circle,inner sep=0pt,minimum size=2.50000000000000pt},>=latex,line join=bevel,]
\node (0) at (7bp,7.5bp) [n] {};
\end{tikzpicture}}
}
,G\hspace{-2pt}\left.
\raisebox{0pt}[3pt][3pt]{}\right)\hspace{-2pt}$

\bigskip

$p_{1}(G)=\hom\hspace{-2pt}\left(\raisebox{0pt}[4pt][4pt]{}\right.
\hspace{-3pt}\ensuremath{\tikzsetfigurename{p1-}
\frac{1}{2}\raisebox{0.5ex}{\tikzsetnextfilename{p1-0}\begin{tikzpicture}[scale=0.150000000000000,baseline=(current bounding box.center),n/.style={fill,shape=circle,inner sep=0pt,minimum size=2.50000000000000pt},>=latex,line join=bevel,]
\node (1) at (78.21bp,7.5bp) [n] {};
  \node (0) at (7bp,18.139bp) [n] {};
  \draw [] (0) ..controls (26.525bp,15.222bp) and (58.714bp,10.413bp)  .. (1);
\end{tikzpicture}}
},G\hspace{-2pt}\left.
\raisebox{0pt}[4pt][4pt]{}\right)\hspace{-2pt}$

\bigskip

$p_{2}(G)=\hom\hspace{-2pt}\left(\raisebox{0pt}[4pt][4pt]{}\right.
\hspace{-3pt}\ensuremath{\tikzsetfigurename{p2-}
\frac{1}{2}\raisebox{0.5ex}{\tikzsetnextfilename{p2-0}\begin{tikzpicture}[scale=0.150000000000000,baseline=(current bounding box.center),n/.style={fill,shape=circle,inner sep=0pt,minimum size=2.50000000000000pt},>=latex,line join=bevel,]
\node (1) at (78.21bp,7.5bp) [n] {};
  \node (0) at (7bp,18.139bp) [n] {};
  \draw [] (0) ..controls (26.525bp,15.222bp) and (58.714bp,10.413bp)  .. (1);
\end{tikzpicture}} + \frac{1}{3}\raisebox{0.5ex}{\tikzsetnextfilename{p2-1}\begin{tikzpicture}[scale=0.150000000000000,baseline=(current bounding box.center),n/.style={fill,shape=circle,inner sep=0pt,minimum size=2.50000000000000pt},>=latex,line join=bevel,]
\node (1) at (79.206bp,69.908bp) [n] {};
  \node (0) at (7bp,69.434bp) [n] {};
  \node (2) at (43.592bp,7.5bp) [n] {};
  \draw [] (0) ..controls (26.798bp,69.564bp) and (59.438bp,69.778bp)  .. (1);
  \draw [] (0) ..controls (17.922bp,50.948bp) and (32.528bp,26.227bp)  .. (2);
  \draw [] (1) ..controls (68.576bp,51.28bp) and (54.361bp,26.37bp)  .. (2);
\end{tikzpicture}}
},G\hspace{-2pt}\left.
\raisebox{0pt}[4pt][4pt]{}\right)\hspace{-2pt}$

\bigskip

$p_{3}(G)=\hom\hspace{-2pt}\left(\raisebox{0pt}[6pt][6pt]{}\right.
\hspace{-4pt}\ensuremath{\tikzsetfigurename{p3-}
\frac{1}{8}\raisebox{0.5ex}{\tikzsetnextfilename{p3-0}\begin{tikzpicture}[scale=0.150000000000000,baseline=(current bounding box.center),n/.style={fill,shape=circle,inner sep=0pt,minimum size=2.50000000000000pt},>=latex,line join=bevel,]
\node (1) at (78.21bp,7.5bp) [n] {};
  \node (0) at (7bp,18.139bp) [n] {};
  \draw [] (0) ..controls (26.525bp,15.222bp) and (58.714bp,10.413bp)  .. (1);
\end{tikzpicture}} + \frac{3}{4}\raisebox{0.5ex}{\tikzsetnextfilename{p3-1}\begin{tikzpicture}[scale=0.150000000000000,baseline=(current bounding box.center),n/.style={fill,shape=circle,inner sep=0pt,minimum size=2.50000000000000pt},>=latex,line join=bevel,]
\node (1) at (149.73bp,20.115bp) [n] {};
  \node (0) at (7bp,17.213bp) [n] {};
  \node (2) at (78.597bp,7.5bp) [n] {};
  \draw [] (1) ..controls (130.23bp,16.656bp) and (98.072bp,10.954bp)  .. (2);
  \draw [] (0) ..controls (26.631bp,14.55bp) and (58.995bp,10.159bp)  .. (2);
\end{tikzpicture}} + \frac{1}{4}\raisebox{0.5ex}{\tikzsetnextfilename{p3-2}\begin{tikzpicture}[scale=0.150000000000000,baseline=(current bounding box.center),n/.style={fill,shape=circle,inner sep=0pt,minimum size=2.50000000000000pt},>=latex,line join=bevel,]
\node (1) at (79.206bp,69.908bp) [n] {};
  \node (0) at (7bp,69.434bp) [n] {};
  \node (2) at (43.592bp,7.5bp) [n] {};
  \draw [] (0) ..controls (26.798bp,69.564bp) and (59.438bp,69.778bp)  .. (1);
  \draw [] (0) ..controls (17.922bp,50.948bp) and (32.528bp,26.227bp)  .. (2);
  \draw [] (1) ..controls (68.576bp,51.28bp) and (54.361bp,26.37bp)  .. (2);
\end{tikzpicture}} - \frac{3}{8}\raisebox{0.5ex}{\tikzsetnextfilename{p3-3}\begin{tikzpicture}[scale=0.150000000000000,baseline=(current bounding box.center),n/.style={fill,shape=circle,inner sep=0pt,minimum size=2.50000000000000pt},>=latex,line join=bevel,]
\node (1) at (98.357bp,65.742bp) [n] {};
  \node (0) at (7bp,43.092bp) [n] {};
  \node (3) at (52.478bp,101.75bp) [n] {};
  \node (2) at (53.166bp,7.5bp) [n] {};
  \draw [] (0) ..controls (30.049bp,48.807bp) and (75.732bp,60.133bp)  .. (1);
  \draw [] (0) ..controls (21.834bp,31.656bp) and (38.099bp,19.116bp)  .. (2);
  \draw [] (1) ..controls (83.114bp,77.705bp) and (67.187bp,90.206bp)  .. (3);
  \draw [] (2) ..controls (52.989bp,31.668bp) and (52.653bp,77.711bp)  .. (3);
\end{tikzpicture}} + \frac{3}{4}\raisebox{0.5ex}{\tikzsetnextfilename{p3-4}\begin{tikzpicture}[scale=0.150000000000000,baseline=(current bounding box.center),n/.style={fill,shape=circle,inner sep=0pt,minimum size=2.50000000000000pt},>=latex,line join=bevel,]
\node (1) at (135.52bp,52.44bp) [n] {};
  \node (0) at (7bp,31.674bp) [n] {};
  \node (3) at (65.649bp,76.768bp) [n] {};
  \node (2) at (76.843bp,7.5bp) [n] {};
  \draw [] (0) ..controls (24.624bp,45.224bp) and (48.41bp,63.513bp)  .. (3);
  \draw [] (1) ..controls (117.89bp,38.936bp) and (94.09bp,20.71bp)  .. (2);
  \draw [] (0) ..controls (26.818bp,24.814bp) and (57.43bp,14.219bp)  .. (2);
  \draw [] (1) ..controls (115.69bp,59.343bp) and (85.07bp,70.006bp)  .. (3);
  \draw [] (2) ..controls (73.625bp,27.415bp) and (68.832bp,57.074bp)  .. (3);
\end{tikzpicture}} - \frac{1}{8}\raisebox{0.5ex}{\tikzsetnextfilename{p3-5}\begin{tikzpicture}[scale=0.150000000000000,baseline=(current bounding box.center),n/.style={fill,shape=circle,inner sep=0pt,minimum size=2.50000000000000pt},>=latex,line join=bevel,]
\node (1) at (92.936bp,56.852bp) [n] {};
  \node (0) at (7bp,44.132bp) [n] {};
  \node (3) at (43.555bp,93.513bp) [n] {};
  \node (2) at (56.291bp,7.5bp) [n] {};
  \draw [] (0) ..controls (19.067bp,60.432bp) and (31.554bp,77.301bp)  .. (3);
  \draw [] (2) ..controls (52.923bp,30.249bp) and (46.985bp,70.35bp)  .. (3);
  \draw [] (0) ..controls (22.504bp,32.609bp) and (40.707bp,19.082bp)  .. (2);
  \draw [] (1) ..controls (77.404bp,68.384bp) and (59.168bp,81.922bp)  .. (3);
  \draw [] (0) ..controls (29.036bp,47.394bp) and (71.018bp,53.608bp)  .. (1);
  \draw [] (1) ..controls (80.84bp,40.561bp) and (68.322bp,23.702bp)  .. (2);
\end{tikzpicture}}
},G\hspace{-2pt}\left.
\raisebox{0pt}[6pt][6pt]{}\right)\hspace{-2pt}$

\bigskip

$p_{4}(G)=\hom\hspace{-2pt}\left(\raisebox{0pt}[8pt][8pt]{}\right.
\hspace{-5pt}\input{includes/p4.inc},G\hspace{-3pt}\left.
\raisebox{0pt}[8pt][8pt]{}\right)\hspace{-2pt}$

\bigskip

$p_{5}(G)=\hom\hspace{-2pt}\left(\raisebox{0pt}[12pt][12pt]{}\right.
\hspace{-6pt}\input{includes/p5.inc},G\hspace{-3pt}\left.
\raisebox{0pt}[12pt][12pt]{}\right)\hspace{-2pt}$

\pagebreak

\end{document}